\documentclass[12pt]{amsart}
\usepackage{amssymb,amsfonts,amsthm,breqn, subfigure,setspace,dcolumn}
\usepackage{cite}
\usepackage{graphicx}
\usepackage{bm}
\theoremstyle{definition}

\theoremstyle{remark}

\numberwithin{equation}{section}
\begin{document}
\title{A fast exact simulation method for a class of Markov jump processes}
\author[Yao Li]{Yao Li}
\address{Yao Li: Department of Mathematics and Statistics, University
  of Massachusetts Amherst, Amherst, MA, 01003, USA} \email{yaoli@math.umass.edu}
\author[Lili Hu]{Lili Hu}
\address{Lili Hu: Georgia Institute of Technology, School of
  Mathematics, 686 Cherry Street, Atlanta, GA, 30332, USA} \email{lilyhu86@gmail.com}
\subjclass[2010]{Primary: 60J22, 65C05, 65C40}

\begin{abstract}
A new method of the stochastic simulation algorithm (SSA), named the
Hashing-Leaping method (HLM), for exact simulations of a class of Markov
jump processes, is presented in this paper. The HLM has
a  conditional constant
computational cost per event, which is independent of the number of
exponential clocks in the Markov process. The main idea of the
HLM is to repeatedly implement a hash-table-like
bucket sort algorithm for all times of occurrence covered by a time
step with length $\tau$. This paper serves as
an introduction to this new  SSA method. We introduce the method,
demonstrate its implementation, analyze its properties, and compare its performance with
 three other commonly used SSA methods in 
  four examples. Our performance tests and CPU operation statistics
show  certain advantage of the  HLM for large scale
problems. 
\end{abstract}


\maketitle

\section{Introduction}

Since the late 1960s, much effort has been devoted to the simulation of Markov
jump processes on high dimensional state spaces. Most of these Markov
jump processes arise in two classes of problems. The first is
usually called chemical reaction networks, which model a
  fixed number of chemical
reactions under dilute, well-mixed conditions. It is well accepted that
when the number of  molecules is small, due to stochastic effects, a deterministic differential
equation fails to model real-world chemical reactions accurately. Therefore, numerous chemical reaction systems
within biological cells, such as
gene networks, regulatory networks, and signaling  pathway networks,
are modeled by Markov
jump processes. The second type of problems are related to the kinetic Monte Carlo
(KMC) method \cite{voter2007introduction}, which essentially covers all
stochastic evolution models that proceed as a sequence of infrequent
transitions at heterogeneous, state-dependent exponential random
times. The KMC was first introduced
to simulate radiation damage \cite{beeler1966displacement}. Today,
it is used to generate stochastic trajectories appearing in
surface/crystal growth, chemical/physical vapor deposition,
 vacancy diffusion,
communication networks, and factory scheduling
\cite{voter2007introduction, petri1962kommunikation,
  kotrla1996numerical, abraham1970computer}.

Markov jump processes coming from both the KMC and chemical reaction networks
have some common features. They are all driven by finitely many
independent exponential clocks. The state of the Markov process is
updated when a clock rings, called an ``event''. The rate of each clock depends on the
current state of the process. Therefore, the update that follows an event may change the rate of other
exponential clocks. In most applications, the update is a simple
transformation, under which the rates of most clocks remain unchanged.
A variety
of Markov jump processes in statistical physics, such as the kinetic
Ising model, the simple inclusion process (SIP),
the simple exclusion process (SEP), and many of their variants like
the TASEP, ASEP, etc. , can also be
categorized into this family.

The scale of these Markov jump processes can be very large. For
example, some chemical reaction networks have thousands of
reactions, while the scale of some reaction-diffusion systems can be as
large as several million. Therefore, it is important to design fast
algorithms for those large scale problems. As a Markov jump process
proceeds sequentially at a series of state-dependent random times, the
fundamental rule of an exact simulation is always to
  identify the time and the index of the next
occurring event, to update the state accordingly,
 and then to advance the time. Algorithms  following this
strategy are called Stochastic Simulation Algorithms (SSA).
Early methods of the SSA like Gillespie's
direct method (DM) \cite{gillespie1977exact}, the first reaction
method (FRM) \cite{gillespie1976general}, and the BKL algorithm
\cite{bortz1975new} for the KMC rely on a linear search of times of
occurrence of events. More sophisticated methods, like the next reaction method (NRM) and
the composition-rejection method (CRM), use a binary heap or other mechanisms to
sample the next occurring event \cite{gibson2000efficient,
  slepoy2008constant}. The performance of a stochastic simulation algorithm is usually measured by computational cost per event. Let
$M$ be the number of exponential clocks in the Markov jump process. Then the computational cost
per event is $O(M)$ for early methods like the DM, the FRM, and the BKL,
$O(\log M)$ for more recent
methods like the NRM, and conditional $O(1)$ for the CRM.  Besides various methods of
SSA, there are also approximate algorithms such as the tau-leaping algorithm and its numerous variants
\cite{cao2006efficient, cao2008slow, gillespie2001approximate}.

A large family of enhanced  methods of the SSA are also developed 
for more specific problems in different applications\cite{cao2004efficient, cao2005multiscale,
  ramaswamy2009new, ramaswamy2010partial, weinan2005nested,
 yang2008kinetic, danos2007scalable}. Important methods that are worth
to mention include: the multi-scale stochastic simulation algorithm (MSSA) for chemical
reaction networks with multiple time scales \cite{cao2005multiscale, weinan2005nested}, the optimized direct method (ODM) for exact simulations of chemical
reaction networks with very heterogeneous reaction rates
\cite{cao2004efficient}, and the next subvolume method (NSM) for stochastic
reaction-diffusion systems \cite{elf2004spontaneous}.

The aim of the present paper is to introduce the Hashing-Leaping
 method (HLM), which is a novel method of the SSA with
conditional $O(1)$ computational complexity per event. Motivated by
the bucket sort algorithm, we repeatedly leap
forward the time by a constant $\tau$, then use a hash-table-like algorithm
to distribute random times covered by the leaping step into $Q$
buckets. Each bucket corresponds to a period of time with length
$\tau/Q$. Under some general assumptions about the Markov process, the
average number of events in each bucket is $O(1)$ for suitable $Q$ and
$\tau$. Then we sequentially
update all events in each bucket until the next leaping step. It is
not difficult to check that the average computational cost per event is $O(1)$ when $\tau \sim O(1)$ and $Q \sim O(M)$. This is further
confirmed by our numerical simulations.

The performance of the HLM is tested in four numerical
examples and compared with that of the DM, the NRM, and the
CRM. The number of
clocks $M$ in the first three examples ranges from tens to
millions. The last example is a chemical oscillator with five
reactions called the ``Oregonator''. Numerical
simulation results show significant advantage of the HLM over the
other tested SSA methods when $M$ is large. For small scale problems,
the HLM remains competitive and is significantly faster than the CRM,
which is the only other existing conditional $O(1)$ method to the best of
our knowledge.

As a novel method of the SSA, there are many extensions that are
beyond the scope of the present paper. Promising extensions include
methods to adjust parameters during the simulation, the
parallelization, extensions to multi-scale problems, and various
applications of the HLM. These topics will be included in our
subsequent works.

The organization of this paper is as follows. In Section 2, we will
describe the family of Markov jump processes that we are interested
in, which covers Markov processes derived from chemical reaction networks
and the KMC. Section 3 presents a short review of
mainstream SSA methods. The HLM is
introduced and analyzed in Section 4. Section 5 focuses
  on performance tests of the HLM on various models. Section 6 is the
conclusion.

\section{Description of Markov jump process model}

We first give a generic description of the Markov jump process to be
studied in the present paper. Let $X_{t} = \{x^{t}_{i}\}_{i = 1}^{N}$ be a Markov jump process on
$\mathbb{R}^{N}$ that is determined by
$M$ random times, which are generated by mutually independent
exponential clocks. The rates of those clocks are state-dependent,
denoted by $R_{1}(X_{t}),
\cdots, R_{M}(X_{t})$, respectively, where $R_{i}: \mathbb{R}^{N} \rightarrow \mathbb{R}^{+}$ are called {\it rate
  functions}. Throughout this paper,  all rate functions are assumed to be time independent. An update transformation $T_{i}^{\omega_{i}}: \mathbb{R}^{N}
\rightarrow \mathbb{R}^{N}$ is associated with each exponential clock,
where $\omega_{i} \in \Omega_{i}$ is a random parameter whose probability
measure is $p_{i}( \mathrm{d} \omega_{i})$, and $\Omega_{i}$ is the sample
space. When the $i$-th clock rings, called an ``event'',
$X_{t}$ is updated by the random transformation $T_{i}^{\omega_{i}}$. During
the update, a random parameter $\omega_{i}$  is sampled from the
probability measure $p_{i}$, independent of everything else. After the update, $X_{t}$ jumps to a new
state $X_{t^{+}} = T_{i}^{\omega_{i}}(X_{t})$. Throughout
the present paper, unless specified otherwise, $M$ means the number of
exponential clocks, which is  said to be the {\it scale} of the Markov jump
process.

The {\it dependency graph} of $X_{t}$ is a directed graph $G = (V, E)$ with $M$
vertices representing $M$ exponential clocks. $\{ i, j \}$ forms a
directed edge if and only if $R_{j} \neq R_{j} \circ T_{i}$. In other
words, $(i, j)$ is an edge if and only if the $i$-th clock affects the
$j$-th clock. 

It is easy to check that Markov jump processes  arising in a very large
family of applications, including statistical mechanics,
chemical reaction networks and  the KMC, fit the description of
$X_{t}$. In those applications, $N$ and $M$ can
be very large numbers, but $R_{i}$ only depends on a limited number of
coordinates of $X_{t}$, and $T_{i}$ is the identity transformation on all but
a limited number of coordinates. Therefore, the dependency graph is
usually sparse,  which means the maximum out degree of $G$
  is independent of $M$. For example, for a
stochastic chemical reaction network with $N$ chemical species and $M$
reactions, we have $X_{t} = \{ x^{t}_{i} \}_{i =
  1}^{N}$, where $x^{t}_{i}$ represents the number of molecules of the $i$-th
reactant species. The rates of reactions, denoted by $\{R_{1}, \cdots,
R_{M} \}$, are determined by the population of reactant
species. The $j$-th transformation is $T_{j}( X_{t}) = X_{t} +
\mathbf{v}_{j}$, where the vector
$\mathbf{v}_{j} \in \mathbb{R}^{N}$ is a sparse vector with all zero
entries except those corresponding to reactant species
  that get changed in the $j$-th reaction.

\section{A short review of existing SSA methods}
\subsection{Direct method  and first reaction method}
As explained in the introduction, when simulating $X_{t}$, it is
important to note that those $M$ clocks are mutually 
independent on a time interval only if all rates $R_{i}$ remain unchanged. When one clock
rings, the corresponding update transformation $T_{i}$ will change the
rates of other clocks. With a small but positive probability, this will
lead to a ``chain reaction'' of events and change the rates of all
clocks in a short time frame. Therefore, to simulate $X_{t}$,
 the strategy of the SSA is  to always
identify the next event.  

The first two popular  methods of the SSA for Markov jump processes
like $X_{t}$ were introduced by Gillespie \cite{gillespie1976general, gillespie1977exact}, called the direct method (DM) and the first reaction
method (FRM), respectively. In the simulation of the
  kinetic Ising model, the DM is also called the BKL algorithm, which was developed
independently \cite{bortz1975new}. In the DM, two random variables are generated
to sample each event. The first random variable determines the time of
occurrence of the next event, and the second one is used to sample the
index of the next event together with a linear search. In the FRM, a linear search is used to
sample both the time of occurrence and the index of the next event
at the same time. Due to the linear search, the computational costs of
both methods are $O(M)$ per event.  These methods can be
summarized as follow. 

{\bf Direct Method}
{\it 
\begin{enumerate}
  \item[{\bf 1: }] Initialize $R_{i}$s. Initialize
      $X_{0}$.  Let $R_{sum} = \sum_{i = 1}^{M} R_{i}$
    and $t = 0$. 
\item[{\bf 2: }] Generate an exponential random variable with rate $R_{sum}$,
  denoted by $\Delta t$
\item[{\bf 3: }] Let $t = t + \Delta t$. Generate a uniform random variable on
  $[0, 1]$, denoted by $u$
\item[{\bf 4: }] Find the minimum $l$ such that $u R_{sum} < \sum_{i = 1}^{l}
  R_{i}$
\item[{\bf 5: }]  Update the state according to
    $T_{l}$: $X_{t^{+}} = T_{l}(X_{t})$
\item[{\bf 6: }] Recalculate all rate functions and $R_{sum}$
\item[{\bf 7: }] Return to {\bf 2} or quit
\end{enumerate}
}
{\bf First Reaction Method}
{\it
\begin{enumerate}
  \item[{\bf 1: }] Initialize $R_{i}$s.  Initialize
      $X_{0}$. Let $t = 0$
\item[{\bf 2: }] Generate $M$ exponential random variables $t_{1}, \cdots, t_{M}$
  with rates $R_{1}, \cdots, R_{M}$, respectively.
\item[{\bf 3: }] Use linear search to find the minimum,
  denoted by $t_{l}$. Let $t = t + t_{l}$
\item[{\bf 4: }]  Update the state according to
    $T_{l}$: $X_{t^{+}} = T_{l}(X_{t})$
\item[{\bf 5: }] Recalculate all rate functions
\item[{\bf 6: }] Return to {\bf 2} or quit
\end{enumerate}
}
It is a standard exercise to show that these two methods are
equivalent. The DM has many variants, such as the
optimized direct method (ODM) introduced by Cao et al \cite{cao2004efficient}. 

\subsection{Next reaction method}

The FRM was significantly optimized by
Gibson and Bruck \cite{gibson2000efficient}, called the next reaction method (NRM). Main
improvements of the NRM include 
\begin{itemize}
  \item  Introducing the concept of the dependency graph. Only rate
    functions affected by an event will be updated
\item Reusing times of occurrence of events without regenerating random
  variables
\item  Using a minimum binary heap to reduce the search time to $O(1)$ and the update time
  to $O(\log M)$
\end{itemize}

The NRM reduces the average computational cost per event to $O(\log
M)$. Currently, this method is widely used in
stochastic simulation packages and commercial software. We
 will also compare the performance of our new method with
 that of the NRM.

The NRM can be summarized as follows
{\it
\begin{enumerate}
  \item[{\bf 1: }] Initialization:
\begin{itemize}
  \item[(a)] Initialize $R_{i}$s.  Initialize $X_{0}$.
\item[(b)] Construct the dependency graph $G$
\item[(c)] Generate $M$ exponential random variables $t_{1}, \cdots, t_{M}$
  with rates $R_{1}, \cdots, R_{M}$, respectively
\item[(d)] Store times of occurrence into a minimum binary heap
\end{itemize}
\item[{\bf 2: }] Find the event on the top of the minimum binary heap, denoted by $t_{l}$
\item[{\bf  3: }]  Update the state according to
    $T_{l}$: $X_{t^{+}} = T_{l}(X_{t})$
\item[{\bf  4: }] Follow the dependency graph to update all affected rate
  function $R_{i}$s
 \item[{\bf  5: }] Update times of occurrence of all affected
   events as
\begin{equation}
  \label{updatetime}
 t^{new}_{i} = (t^{old}_{i} - t_{l})\cdot
  \frac{R^{old}_{i}}{R^{new}_{i}} + t_{l}  \,,
\end{equation}
and maintain the binary heap
\item[{\bf 6: }] Advance $t_{l}$ by an exponential random
  number with rate $R_{l}$ and maintain the binary
  heap
\item[{\bf 7: }] Return to {\bf 2} or quit.
\end{enumerate}
}

Equation \eqref{updatetime} is also adopted by the CRM \cite{slepoy2008constant} and the HLM
introduced in this paper. As will be explained in Section 4.2,
equation \eqref{updatetime} comes
from the property of the exponential distribution. We remark that this
transformation formula only applies to time independent rate
functions. Equation \eqref{updatetime} will be different if some
$R_{i}$s are time varying, see \cite{gibson2000efficient} for the
detail. 

\subsection{Composition-Rejection method}

The idea of the composition-rejection method (CRM)
\cite{slepoy2008constant} comes from a simple probabilistic fact
 implicit in Gillespie's direct method. Suppose we have $M$ 
mutually independent exponential
random variables $Y_{1}, \cdots, Y_{M}$ with 
rates $R_{1}, \cdots, R_{M}$, respectively. Let $Y_{min}$ be the
minimum of the $M$ random variables and $R_{sum}$ be the sum of the
$M$ rates. Then $Y_{min}$ is an exponential random variable with
rate $R_{sum}$. In addition,
$$
  \mathbb{P}[ Y_{min} = Y_{i}] = \frac{R_{i}}{R_{sum}} \,.
$$

Therefore, sampling $Y_{i} = Y_{min}$ is equivalent to sampling a weighted
distribution over $\{1, \cdots, M \}$ with weights $R_{1}/R_{sum}, \cdots, R_{M}/R_{sum}$,
respectively. Instead of the linear search used in the DM or the FRM, this
sampling can also be done by the following rejection-based method. Let
$R_{max} = \max \{R_{1}, \cdots, R_{M} \}$. Two random numbers $Z_{1}$
and $Z_{2}$ are repeatedly generated until $Y_{min}$ is selected, where $Z_{1}$ is
uniformly distributed on $[0, R_{max}]$, and $Z_{2}$ is an integer that
is uniformly distributed over $\{1, \cdots, M \}$. The rule of the rejection is
that: if $Z_{1} > R_{Z_{2}}$, then the pair $(Z_{1}, Z_{2})$ is
rejected. Otherwise it is accepted and we have $Y_{min} = Y_{Z_{2}}$. 

The rejection-based sampling method has constant computational
complexity independent of $M$. However, the constant can be very large
if $R_{max}$ is much larger than all the other $R_{i}$s. This is
partially solved by the CRM. The CRM relies on
the maximal and minimal rates $R_{min}$ and $R_{max}$. To reduce the
expected number of rejections, $R_{1}, \cdots,
R_{M}$ are distributed into $g \geq \left \lceil \log_{2} (R_{max}/R_{min})
  \right \rceil$ groups. The first group
contains rates ranging from $R_{min}$ to $2 R_{min}$, the second group
from $2 R_{min} $ to $4R_{min}$, and so on. According to \cite{slepoy2008constant},
$g \leq 30$ is sufficient for most applications. Sums of rates in each
groups are calculated, denoted by $p_{1}, \cdots, p_{g}$, and
stored. The CRM uses Gillespie's direct method to  sample the
group (composition), and uses the rejection-based sampling technique
introduced above to select the event (rejection) within the sampled
group. After updating an event, $R_{sum}$, $R_{max}$, $R_{min}$ and all affected
groups will be maintained.  

The CRM reduces the computational cost
per event to conditional $O(1)$ when $M$ is sufficiently
  large. See the performance test in Section 5.2. The requirement of
the $O(1)$ complexity is (i) the dependency graph $G$ is sparse and (ii) the ratio $R_{max}/R_{min}$ is $O(1)$.

The  CRM can be summarized as follows:
{\it
\begin{enumerate}
  \item[{\bf 1: }] Initialization:
\begin{itemize}
  \item[(a)] Initialize $R_{i}$s. Initialize $X_{0}$. Let
    $t = 0$
\item[(b)] Construct the dependency graph $G$
\item[(c)] Distribute $M$ clock rates into $g$ groups and compute  $p_{1},
  \cdots, p_{g}$, the sum of
  rates in each groups
\item[(d)] Calculate $R_{sum} = \sum_{i = 1}^{M} R_{i}$
\end{itemize}
\item[{\bf  2: }] Generate an exponential random variable with rate $R_{sum}$,
  denoted by $\Delta t$
\item[{\bf  3: }] Let $t = t + \Delta t$. Generate a uniform random variable on
  $[0, 1]$, denoted by $u$
\item[{\bf  4: }] Find the minimum $k$ such that $u R_{sum} < \sum_{i = 1}^{k}
  p_{i}$
\item[{\bf  5: }] Use the rejection-based sampling method to select an event $R_{l}$ from
  group $p_{k}$
\item[{\bf  6: }]   Update the state according to
    $T_{l}$: $X_{t^{+}} = T_{l}(X_{t})$
\item[{\bf  7: }] Follow the dependency graph to update all affected rate
  function $R_{i}$s
 \item[{\bf  8: }] Update times of occurrence of all affected
   events as \eqref{updatetime}
\item[{\bf  9: }] Maintain the groups, update
  $R_{sum}$ , $R_{max}$, and $R_{min}$ and affected $p_{i}$s
\item[{\bf  10: }] Return to {\bf 2} or quit.
\end{enumerate}

}

\section{Hashing-Leaping method (HLM)}

\subsection{Introduction to the  method}

In this section, we introduce a  conditional
  $O(1)$ per-event  SSA method for the exact
simulation of $X_{t}$. As reviewed in the previous section, the main bottleneck of
simulating $X_{t}$ is sampling the ``next event''. Instead of
a linear search or a heap sort, the HLM is motivated by the bucket
sort algorithm, which is a linear complexity 
sorting algorithm in most practical settings. To simulate $X_{t}$, two
parameters $\tau$ and $Q$ are chosen by either observing rate
function $R_{i}$s or  performing a smaller scale simulation, where $\tau > 0$ is the
step size and the positive integer $Q$ is the number of buckets. 

The HLM runs in the following way:  Same as in the
  NRM, times of occurrence of
events associated with $M$ clocks, denoted by $t_{1}, \cdots, t_{M}$, are stored and
maintained. The algorithm makes a major update in the beginning of
every time step with length $\tau$, called a {\it bucket
  redistribution}. In the $n$-th bucket redistribution, $t_{1},
\cdots, t_{M}$ are distributed into $Q+1$ buckets, denoted by $B_{1},
\cdots, B_{Q}, B_{L}$, that represent time intervals
\begin{eqnarray*}
B_{1} & = &  [ (n-1)\tau  \, , \, (n-1)\tau + \tau/Q )\\
B_{2} & = & [  (n-1)\tau + \tau/Q \, , \, (n-1)\tau + 2\tau/Q )\\
&\vdots&\\
B_{Q} & = & [ (n-1)\tau + (Q-1)\tau/Q \, , \,
n\tau ) \\
B_{L} & = & [ n \tau \,, \, + \infty ) \, ,
\end{eqnarray*}

respectively. Then we start from the first
non-empty bucket $B_{n_{1}}$ to find the minimum time of occurrence, say $t_{l}$, by a linear search, and make update according to $T_{l}$. During the update, the following two operations will be
carried out:  (i) The dependency graph is followed to
update the new clock rates of all affected clocks, as well as the
corresponding times of occurrence. (ii) An exponential random variable with rate $R_{l}$
will then be generated and added to $t_{l}$, which is the time of
occurrence of the next event associated with $R_{l}$. $t_{l}$ will then be placed
into the proper bucket. We repeat this step until
  $B_{n_{1}}$ is emptied.

Then we move to the next nonempty bucket and carry out the same series
of 
operations. This procedure
continues until all buckets $B_{1}, \cdots, B_{Q}$ are emptied. At that
time, times of occurrence of all events are stored in $B_{L}$,
 which will be used to perform the next
bucket redistribution. We call it the Hashing-Leaping method because the bucket
redistribution step resembles the Hash algorithm, while the whole
algorithm can be seen as an exact version of the tau-leaping
algorithm.

The HLM should be implemented with the proper data structure
to improve the efficiency. The simplest way we find is to construct an array of
$N$ structs, named the {\bf TimeArray}. Each 
struct, of the type {\bf ST}, has three elements: a floating
point number that
indicates the time of occurrence of the event associated
with the
exponential clock, and  two {\bf ST} pointers pointing to
  its left and right neighbors, respectively. In addition we need an array of $Q+1$
{\bf ST} pointers, called the {\bf BucketArray}, that represents the
heads of the $Q+1$
buckets. A floating point number array {\bf RateArray} is also needed to store
the rate of each clock. 

\begin{figure}[h]
\centerline{\includegraphics[width = \linewidth]{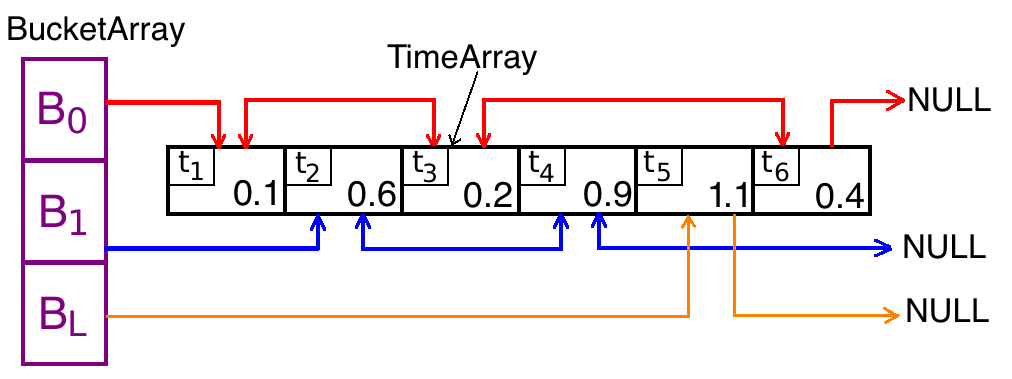}}
\caption{TimeArray and BucketArray for $M = 6$, $Q = 2$ and $\tau =
  1$. The TimeArray stores times of occurrence of all $6$
  clocks, labeled $t_{1}$ to $t_{6}$. Three buckets stored in the BucketArray, labeled $B_{0}$, $B_{1}$, and
  $B_{L}$, contain times of occurrence between $0$ and $0.5$, between $0.5$
and $1.0$, and greater than $1.0$, respectively. Arrows represent
the direction of pointers}
\label{algorithm} 
\end{figure}

After a bucket redistribution, each bucket is formed by a doubly linked list whose
head is pointed to by an element in the {\bf BucketArray}, as shown in Figure
\ref{algorithm}. When updating each
bucket, a linear search is performed to find the minimum time of occurrence
within this bucket. Every update after an event requires two operations: (i) remove the
corresponding {\bf ST} struct from its old bucket, i.e., set its left
and right pointers  to  NULL and maintain the doubly linked list; and
(ii)  push the struct into the front of the new bucket, i.e., relink its
left and right pointers. To increase the efficiency, if an {\bf ST}
struct remains in the same bucket after an update, only its time of
occurrence will be changed. It is a simple practice to implement the
HLM in C/C\verb!++!.  

We remark that according to our test, it seems to be less efficient to
implement buckets as linear arrays than to implement them as linked lists. Although linear
arrays are more cache friendly than linked lists, we have to
frequently move structs from one bucket to the other instead of just
relinking pointers. In addition, to maintain the linear array data
structure, when removing one element from the bucket, the last element
has to be moved to fill the empty slot. As a result, we observed $\sim 10 \%$ decrease of
the performance when implementing buckets as linear arrays. 

The HLM can be summarized as follows
{\it 
\begin{enumerate}
  \item[{\bf 1: }] Initialization:
\begin{itemize}
\item[(a)] Initialize $R_{i}$s. Initialize $X_{0}$. Let $t = 0$.
  \item[(b)] Generate the dependency graph $G$. 
\item[(c)] Generate times of occurrence $t_{1}, \cdots, t_{M}$.
\item[(d)] Initialize {\bf TimeArray}, {\bf RateArray}, and
    {\bf BucketArray}
\end{itemize}
     
\item[{\bf  2: }] Choose proper parameters $\tau$ and $Q$.
\item[{\bf  3: }] Distribute times of occurrence to corresponding
  buckets {\bf BucketArray[0]} $\sim$ {\bf BucketArray[Q]}.
\item[{\bf  4: }] For $i = 0$ to $Q - 1$

\hskip -0.5cm While {\bf BucketArray[i]} is nonempty: 

\begin{itemize}
  \item[(a)] Find the least time of occurrence $t_{l}$ within this bucket
\item[(b)] Make update according to $T_{l}$
\item[(c)] Follow the dependency graph $G$ to update all affected
  rate functions in {\bf RateArray}
\item[(d)]  Update time of occurrence of all affected events as \eqref{updatetime}
and move affected elements in {\bf TimeArray} to the
corresponding buckets. 
\item[(e)] Advance $t_{l}$ by an exponential random number with rate $R_{l}$ and move {\bf TimeArray}[l] to the
  corresponding bucket
\end{itemize}

\item[{\bf  5: }] Update the time intervals of each bucket and return to
  {\bf 3}, or quit.
\end{enumerate}
}

\subsection{Analysis of algorithm}

It is important to demonstrate the correctness of the HLM before
further investigations. In fact,
the HLM is mathematically equivalent to the NRM. This can be checked
by using the
following three steps.

\begin{enumerate}
  \item[{\bf 1: }] {\bf Sampling events } Sampling the next time of
    occurrence is the core step of the SSA. In the HLM, the next time
    of occurrence is chosen as the minimum of $M$ random times, which is
    equivalent to the NRM and the FRM.

\item[{\bf 2: }] {\bf Reusing random times.} After initialization, a sample of the time
  of occurrence of each clock is taken from the corresponding
  exponential distribution. If an event occurs at time $t$ without
  changing the rate of this clock, this sample can be reused due to
  the  time-invariant nature of the exponential distribution. More precisely, if $Z$ is an exponentially
  distributed random variable, then for any $t > 0$, the
  ``overshoot'', i.e., $Z |_{Z > t}$, has the same exponential
  distribution. We refer Theorem 1 in \cite{gibson2000efficient} for
  the full mathematical detail.

\item[{\bf 3: }] {\bf Changing clock rates.} The update transformation in an event may
  change the rate of other clocks. Therefore the affected times of
  occurrence need to be updated accordingly. 

  Assume after an event at time $t$, the rate of the $i$-th clock changes from
  $R^{old}_{i}$ to $R^{new}_{i}$. By step {\bf 2}, the ``overshoot''
  $t_{i} - t$ has an exponential distribution with rate
  $R^{old}_{i}$. In addition, it is a simple probabilistic fact that for any
  constant $a > 0$ and any exponential random variable $Z$ with rate
  $r$, $a Z$ is exponentially distributed with rate $r/a$. Hence
  the transformation in equation \eqref{updatetime} maps an
  exponential distribution with rate $R^{old}_{i}$ to an exponential
  distribution with rate $R^{new}_{i}$. For further reference
  regarding the proof, see Theorem 2 in
  \cite{gibson2000efficient}.

\end{enumerate}

\subsection{Analysis of complexity}
\begin{enumerate}
  \item 

{\bf Average computational cost. }

 Assume 
\begin{itemize}
  \item[(a)] There exists a constant $K$ independent of $M$ such that
$$
  \mathbb{E}[\mbox{Number of events occurring on } [t \, , t + \Delta] ] 
  \leq KM\Delta
$$ 
for any $t$ and any $\Delta > 0$; and 
\item[(b)] The dependence graph is sparse such that the maximum out
  degree is $O(1)$.
\end{itemize}

Then if $\tau \sim O(1)$, $O(M)$ events occur in
each time step with length $\tau$. The 
computational cost of a bucket redistribution is $M$. If we choose $Q
\sim O(M)$, then the number of events in each
bucket can be (very) roughly approximated by a Poisson distribution
with $O(1)$ mean. Since the dependency graph is sparse, the
average cost of updating each bucket is $O(1)$. Therefore the total computational cost in one step is $O(M) +
Q \times O(1)= O(M)$. This makes the average complexity of the HLM be
$O(1)$ per event. 

If the dependency graph is not sparse in a way that the average out degree of
vertices is $D$, which is possibly dependent of $M$, then the average
cost of updating each bucket is $O(D)$. This brings the total
computational cost in one step to $O(D M)$ and the average complexity
per event to $O(D)$.

We remark that assumptions (a) and (b) are satisfied by a large class
of models in practice. In particular, if $R_{max} \sim O(1)$ and
$R_{max}/R_{min} \sim O(1)$, then (a) is satisfied.
\item 

{\bf Worst case analysis. } In the worst case, which means all events
are distributed into the same bucket, the HLM is equivalent
to the FRM. Therefore, the computational cost per event
in the worst case is $O(M)$. However, we remark that  in
  most practical cases,  this worst situation
occurs with an extremely low probability. If 
  $R_{max} \sim O(1)$ and
$R_{max}/R_{min} \sim O(1)$, when letting $Q \sim O(M)$ and
$\tau\sim O(1)$, the probability that all $M$ events are placed into a
single bucket is $\sim O(M^{-M})$.   

\end{enumerate}

\subsection{Discussion of issues}
\begin{enumerate}
\item 

{\bf Choice of parameters. } 

Finding optimal parameters for the HLM is a challenging
  job. We make some idealized calculations to shed some light on this
  problem. Assume (a) and (b) in Section 4.3 hold. Then the average computational cost per event, denoted by $\mathcal{C}$, on this time step satisfies
\begin{eqnarray*}
 \mathcal{C} & =& \frac{1}{\# events}  \bigg \{ search \,\, cost +
     update \,\, cost  \\
&& + bucket \,\, iteration \,\, cost + bucket \,\,
       redistribution \,\, cost \bigg \}\\
&=&\frac{O(1)}{\alpha M \tau} \bigg \{ C_{s} Q \cdot \frac{1}{2}\left [ (\frac{\alpha M
      \tau}{Q})^{2} + 2\frac{\alpha M
      \tau}{Q}  \right ]  + C_{u}\alpha M \tau\\
&&  + C_{i} Q + C_{r}M
  \bigg \} \,,
\end{eqnarray*}
where $C_{s}$, $C_{u}$, $C_{i}$, and $C_{r}$ represent the average
cost of searching, updating events, bucket iteration, and bucket
redistribution. Therefore, it is easy to check that for each fixed
$\tau$, the optimal $Q$ satisfies 
$$
  Q_{opt} = \alpha M \tau \cdot \sqrt{\frac{C_{s}}{2 C_{i}}}  \,.
$$
When $Q$ is optimal, we have
$$
  \mathcal{C} = \sqrt{2 C_{s} C_{i}} + C_{s} + C_{u} +
  \frac{C_{r}}{\alpha \tau} \,.
$$
It is reasonable to assume $C_{s}$, $C_{i}$, and $C_{r}$ are
constants. But $C_{u}$ depends on $\tau$ because we do not relink pointers if an update does not
move the affected event to a new bucket. For smaller $\tau$, larger
proportion of events will remain in the bucket $B_{L}$ throughout
the step. Hence a more precise model of
$C_{u}$ is
$$
  C_{u} = C_{u}(\tau) =  C'_{u} + C''_{u}(1 - e^{-\alpha \tau} ) \,,
$$
where $C'_{u}$ and $C''_{u}$ are constants that represent the cost of update without
relinking pointers and the cost of relinking pointers.

Therefore, we conclude that the optimal number of buckets $Q$ is in proportion to
both $\tau$ and $M$, and the optimal value of $\tau$ depends on constants
$C''_{u}$, $C_{r}$, and $\alpha$. 

One important remark is that the empirical performance of the HLM is not
sensitive with respect to small change of parameters. The reason is that
updating an event, which includes calculating new rate functions,
modifying times of occurrence, and relinking pointers, is much more
expensive than the other operations. Hence $C_{u}$ is significantly greater than the other three constants. Similarly,
$C'_{u}$ is greater than $C''_{u}$. According to our test, the
performance of the HLM is stable as long as $Q \sim
O(M)$ and $\tau \sim O(1)$. For example, for the generalized KMP model
introduced in Section 5.1, the CPU time with $\tau = 1, Q
= M$ and that with $\tau = 1, Q = M/10$ have only less than $ 10 \%$
difference. 

When simulating models in Section 5, we choose smaller $\tau$ such that
about half of events are placed into $B_{L}$, then make $Q$
be in proportion to $\tau$ and $M$. We admit that these parameters may
not be optimal. To obtain the optimal
parameters, all constants introduced above should be estimated
empirically. In other words, the algorithm needs to be enhanced such that $\tau$ and $Q$ can be properly adjusted as
the simulation proceeds. We will make this extension in our subsequent works.

  \item {\bf Possible improvements. } 

There are two places to further
    improve the HLM. The first potential
improvement is at the level of implementation. One can replace the
linear search in each bucket with a binary search. Under our assumption, the expectation of
the largest number of events stored in a bucket is $O(\log
M)$. Therefore, a binary search could potentially 
increase the performance of the algorithm. As mentioned before, it is slightly
less efficient to implement buckets as linear arrays than to implement
them as linked lists. However, we
expect some improvement if buckets are implemented as linear arrays, with
binary heaps constructed on them. ( It is standard knowledge that implementing a binary heap as a
linked list is significantly more complicated than implementing that
as a linear array. ) We did not present this
improvement in the present paper because it does not change the $O(1)$
computational complexity per event of the HLM. Plus due to the overhead of
constructing the binary heap, we expect the improvement can only be
observed for very large scale problems.

The second improvement is for Markov jump processes with sparse
dependency graphs. Assuming $R_{max}/R_{min} \sim O(1)$,  when the dependency graph is sparse,
 times of occurrence stored in each bucket will not affect each
other with high probability ($ \sim 1 - O(1/M)$). Therefore, in
most situations we can update times of occurrence in each bucket sequentially from
the head of the list without searching for the minimum. Since in average a bucket only
contains $O(1)$ events, this improvement does not change the computational
complexity of the algorithm either. In fact, we only observed a minor increase
of simulation speed at the cost of much more complicated programming
code. However, this idea can be used in the parallelization of our
algorithm, as explained in (3).

\item {\bf Parallelization.}

Although the fact that events in one bucket are ``almost independent''
does not significantly improve the performance of the HLM when
running on
a single CPU, it can be used to parallelize this algorithm. In
fact, the HLM is very compatible with parallel programming if
the dependency graph is sparse. The idea is to make the
number of buckets $Q \sim O(M/p)$, where $p$ is the number of CPUs, and divide one bucket into
$p$ sub-buckets. During the bucket redistribution, we can evenly
distribute events into sub-buckets that are maintained by different CPUs. Suppose the dependency graph is
sparse, then events placed in each bucket will not affect each other
with high probability. Therefore in most of the running time, different CPUs can update events in their own
sub-buckets independently.  

The principle of the parallelization is straightforward. However, there
are lots of details remaining to be studied. With a low but strictly positive probability, an event in one
bucket can affect events in the same bucket (but probably different
sub-bucket), which will cause intensive communication between
CPUs. Therefore, it is crucial to design algorithms that can 
identify and deal with these small probability events with the lowest overall
overhead. We will study the parallel HLM and its applications intensively in our subsequent
works.

\end{enumerate}

\section{Numerical Examples}

\subsection{Introduction to models}

We chose the following four models to test the performance of the
 HLM. 

{\bf (i). Generalized KMP model.}
There are numerous stochastic processes in the field of statistical
mechanics that fit the generic description of $X_{t}$ in Section 2,
such as the simple inclusion process (SIP), the simple exclusion process
(SEP) and its variants \cite{boldrighini1989computer,
  derrida1993exact, derrida1998exactly}, and many stochastic processes that models the
microscopic heat conduction \cite{eckmann2006nonequilibrium,
  li2014nonequilibrium, bonetto2000fourier}. We chose the following generalized KMP
model as a test case of our simulation. The KMP model is a stochastic
model proposed in \cite{kipnis1982heat} that models the microscopic energy
transport, from which a microscopic derivation of Fourier's law was
carried out. The main feature of the generalized KMP model ( See
\cite{li2013existence, grigo2012mixing} for details ) here is its
energy dependent clock rates, which makes the simulation non-trivial.


The KMP model models the energy transport in a 1D chain of $M$
oscillators coupled with two heat baths whose temperatures are $\mathcal{T}_{L}$ and
$\mathcal{T}_{R}$, respectively. We only take note of energy carried by each
oscillator, denoted by $x_{1}, \cdots, x_{M}$, respectively. The Markov chain
generated by the KMP model can be described as follows. An exponential
clock with rate $R(x_{i}, x_{i+1}), i = 0 \sim M$ is
associated with each pair of adjacent oscillators (let $x_{0} = \mathcal{T}_{L}$
and $x_{M+1} = \mathcal{T}_{R}$). When the clock rings, the energy stored in
the corresponding pair of oscillators is pooled together, repartitioned
randomly, and redistributed back to the oscillators. The energy
redistribution satisfies 
$$
  (x_{i}', x'_{i+1}) = \left ( p(x_{i} + x_{i+1}), (1-p)(x_{i} +
    x_{i+1}) \right ) \,,
$$
where $p$ is a uniform random variable on $(0,1)$ that is independent
of everything else, and $x_{i}'$ is the energy carried by oscillator $i$ after the
update. If clock $0$ (resp. clock $M$) rings, the energy of the first (resp. last)
oscillator exchanges energy with an exponential random variable with
mean $\mathcal{T}_{L}$ (resp. $\mathcal{T}_{R}$) in the same way.  

We let $\mathcal{T}_{L} = 1.0$, $\mathcal{T}_{R} = 2.0$, and $R(x_{i}, x_{i+1}) =
\sqrt{x_{i} + x_{i+1}}$ and simulate the generalized KMP model with varying $M$s.  

Markov jump processes arising in other
statistical mechanics models may not satisfy the
generic description of $X_{t}$ in Section 2. However, they can still be
simulated by the HLM efficiently. One example is the random
halves model proposed in \cite{eckmann2006nonequilibrium} and its
modification in \cite{li2014nonequilibrium}, in
which the number of clocks varies randomly with time instead of being a
constant. When the number of exponential clocks changes by one, the cost of
maintaining  the bucket data structure is $O(1)$. In contrast, in the
NRM, the computational cost of maintaining the binary heap after
adding/removing one exponential clock is $O(\log M)$.

{\bf  (ii). Chemical reaction network} The second example we chose is a
chemical reaction network. Stochastic chemical reaction networks are a
very important class of Markov jump process. It is
well known that  most SSA methods were originally developed for the stochastic simulation of chemical
reactions. To simplify the implementation, we use
the same example in \cite{slepoy2008constant} with minor
modification.

Consider a large chemical reaction network with $M$ reactions. For the
sake of simplicity we only take note of rate functions, or the
propensity. The initial propensities are uniformly distributed over
$[0, \,  2]$. The dependency graph is randomly generated in a way that
each reaction affects $m$ other reactions, where $m$ is a random
integer uniformly distributed over $\{1, \cdots, 30 \}$. After the
occurrence of each
reaction, all affected reactions are updated in a way that the
propensity is replaced by a new random number uniformly distributed on
$[0, \, 2]$. We note that this example is 
different from real-world chemical reaction networks, and is
used only for the purpose of testing algorithms. 

The main difference between our example and the example
presented in \cite{slepoy2008constant} is that we do not bound the
ratio of $R_{max}$ to $R_{min}$. In the example used by
\cite{slepoy2008constant}, the ratio of $R_{max}$ to $R_{min}$ is
bounded by $1 \times 10^{6}$ at the initial state and has minor
changes when the simulation proceeds.

\begin{figure}[h]
\centering
\subfigure{\includegraphics[width = 0.48
  \linewidth]{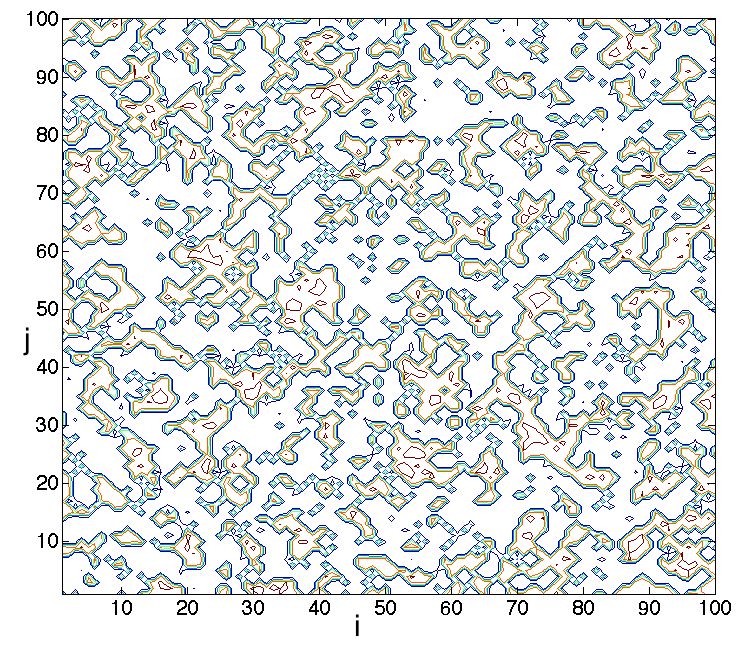}} \subfigure{\includegraphics[width = 0.51 \linewidth]{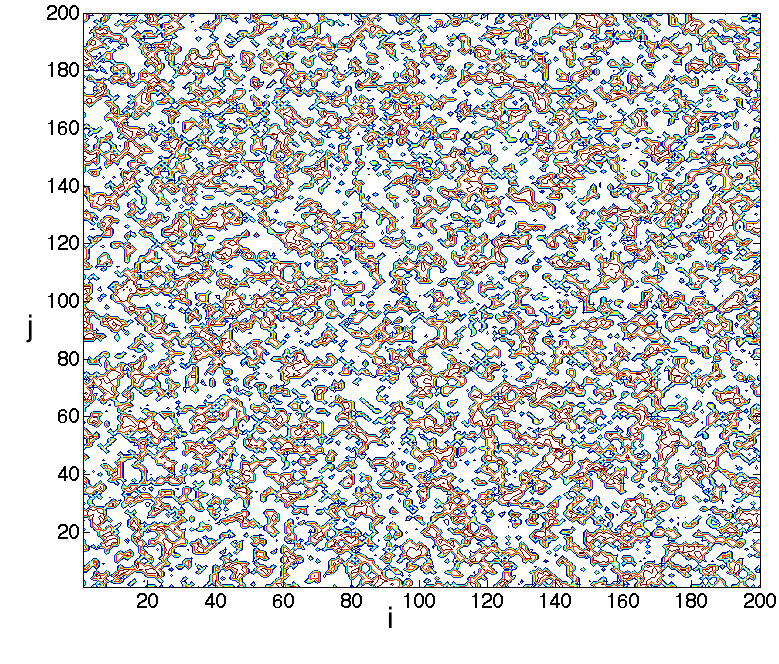}}

\caption{Pattern formation of the Gray-Scott model. Figure 2 (a):  Contour plots of
  numbers of $U$ molecules on $100 \times 100$ subvolumes at $t =
  1500.0$. Figure 2(b): Contour plots of
  numbers of $U$ molecules on $200 \times 200$ subvolumes at $t =
  1500.0$}
\label{pattern}
\end{figure}

{\bf (iii). Reaction-diffusion system}
The third example we chose is the Gray-Scott model, which is a
reaction-diffusion system  in 2D.  Since the purpose of this
  simulation is to compare  the performance of algorithms, we simplify some
details. The Gray-Scott model is well known as its
 pattern formation phenomenon \cite{vigelius2012stochastic}, as seen in Fig \ref{pattern}. The
reaction part of the model involves two species $U$ and $V$ and four
reactions 
\begin{eqnarray*}
U + 2 V & \xrightarrow{K_{1}(\Omega)} & 3 V \\
U & \xrightarrow{K_{f}} & \emptyset\\
V & \xrightarrow{K_{2}} & \emptyset\\
\emptyset & \xrightarrow{K_{f} u_{0}(\Omega)}& U
\end{eqnarray*}
where $\Omega$ is a parameter that indicates the scale of the system
 per subvolume. $\Omega$ is determined by both the molecular
  population level and the edge length of subvolumes, which are
  assumed to be constants throughout this example. Parameters
  $k_{1}(\Omega) = \hat{k}_{1} \Omega^{-2}$ and $u_{0}( \Omega) = 
\hat{u}_{0} \Omega$ depend on $\Omega$. $\hat{k}_{1}, K_{f}, K_{2},
\hat{u}_{0}$ are constants. Molecules are assumed to be well-mixed in each
subvolume. Besides reactions, $U$ and $V$ molecules can
jump to neighbor subvolumes at  certain diffusion rates, denoted by
constants $D_{U}$ and $D_{V}$, respectively.  $D_{U}$ and
  $D_{V}$ also depend on the edge length of
  subvolumes. 

In the simulation, we study domains that are consist of $K
  \times K$ subvolumes for $K$ ranging from $3$ to $1200$.  Numbers of molecules of species $U$ and $V$ in each
subvolume are denoted by $U_{i,j}$ and $V_{i,j}$, respectively, where $i,
j$ ranges from $1$ to $K$. There are six exponential clocks 
associated with each subvolume with rates $U_{i,j}V_{i,j}^{2}K_{1}(
\Omega)$, $U_{i,j} K_{f}$, $V_{i,j} K_{f}$, $K_{f}u_{0}(\Omega)$,
$U_{i,j}D_{U}$, and $V_{i,j} D_{V}$, respectively.  When
  the diffusion clock rings, one corresponding molecule in the subvolume moves to one of its nearest neighbors with
equal probability. A molecule exits from the system if it moves out of
the domain. Hence the scale of this
system is $M = 6K^{2}$. The values of parameters are
taken as $\Omega = 250$, $K_{f} = 0.0055$, $K_{2} = 0.0205$, $D_{V} =
0.002$, $D_{U} = 0.001$, and $\hat{k}_{1} = \hat{u}_{0} =
1$.

{\bf (iv). Chemical oscillator system.}

To examine the performance of the HLM for small scale problems, we
consider the following chemical oscillator system called the
``Oregonator'', which was originally introduced in \cite{field1974oscillations}. The
``Oregonator'' is a chemical reaction system with five reactions
\begin{eqnarray*}
\bar{X}_{1} + Y_{2} & \xrightarrow{c_{1}} & Y_{1} \\
Y_{1} + Y_{2} & \xrightarrow{c_{2}} & Z_{1}\\
\bar{X}_{2} + Y_{1} & \xrightarrow{c_{3}} & 2 Y_{1} + Y_{3}\\
2 Y_{1} &\xrightarrow{c_{4}} & Z_{2}\\
\bar{X}_{3} + Y_{3} & \xrightarrow{c_{5}} & Y_{2} \,,
\end{eqnarray*}
where the molecular population level of species $\bar{X}_{1},
\bar{X}_{2}$, and $\bar{X}_{3}$ are assumed to be constants. The
parameters are taken to be $c_{1} \bar{X}_{1} = 2$, $c_{2} = 0.1$,
$c_{3} \bar{X}_{2} = 104$, $c_{4} = 0.016$, and $c_{2} \bar{X}_{3} =
26$, which are consistent with \cite{gillespie1977exact}. 

\subsection{Performance test of HLM}

{\bf (i). Speed of simulation}

The first set of numerical simulations concern the speed of algorithm. We
implement  each of the four models introduced in the previous subsection using four
different methods: the CRM, the DM, the HLM, and the NRM. All
SSA methods 
are implemented in C, with C\verb!++! I/O for the sake of simplicity of
coding. Implementations of all algorithms are
  optimized to the best of our ability. For example, in the
  implementation of the CRM, instead of dynamic group bounds varying with $R_{min}$,
  pre-assigned group bounds are used to reduce the overhead of maintaining 
groups.   All performance tests are run on a 2012 Macbook Pro
laptop with an Intel Core I7-3615QM CPU and $8$ GB memory.

Instead of stopping after simulating a fixed number of events, we chose to simulate
all examples up to $t = 10$, as simulating Markov processes up to
different times may bias the result. The performance
   of four SSA methods is
measured in seconds per million events and compared. The scales
of the generalized KMP model, the chemical reaction network, and the
reaction-diffusion system ranges from $M = 10 \sim 10^{7}$, $M =
50 \sim 5 \times 10^{6}$, and $M = 54 \sim 8.64 \times 10^{6}$,
respectively. 

The parameters of the HLM are chosen as $\tau = 0.2$, $Q = M/10$ for the generalized
KMP model, $\tau = 0.1$, $Q = M/20$ for the chemical reaction network,
$\tau = 0.5$, $Q = M/2$ for the reaction-diffusion system,  and $\tau =
0.01$, $Q = 5$ for the ``Oregonator''. The parameter of the CRM is chosen
as $g = 30$ for all models.

\begin{figure}[h]
\centering
\includegraphics[width = \textwidth]{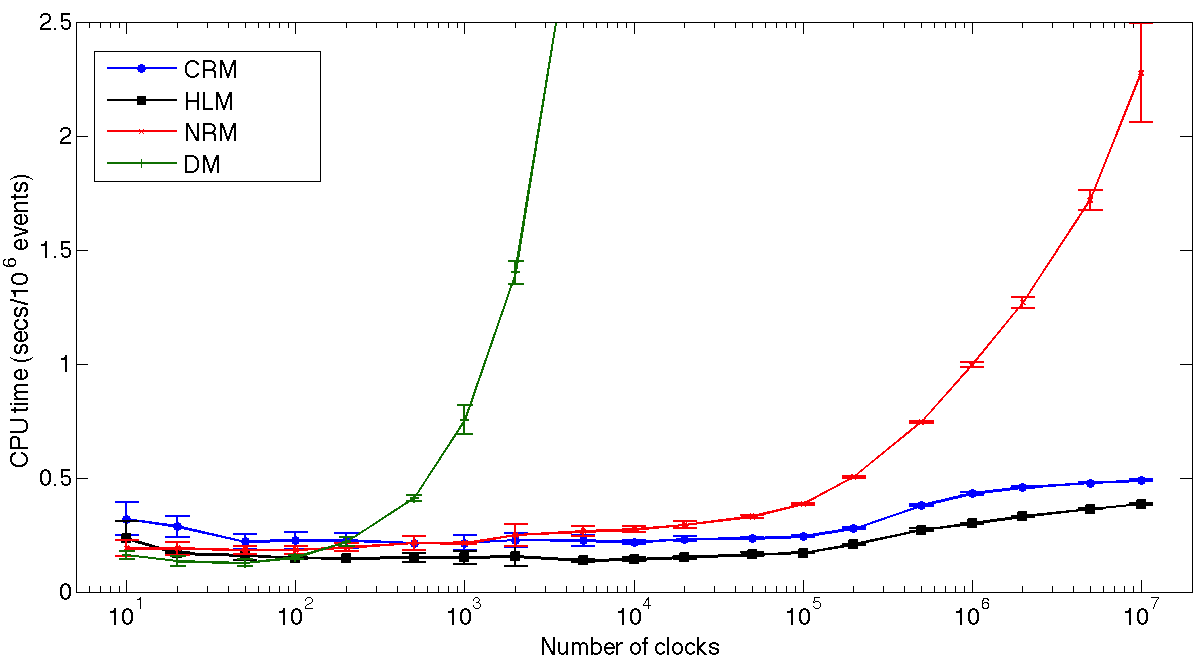}
\caption{CPU times in seconds per $10^{6}$ events vs. number of
  exponential clocks for the generalized KMP model. Blue, black, red, and green
  plots represent the mean CPU times of the CRM, the HLM, the NRM,
  and the DM over $10$ runs, respectively. The error bars indicate
  one standard deviation of the mean.}
\label{KMPspeed}
\end{figure}

\begin{figure}[h]
\centering
\includegraphics[width =  \textwidth]{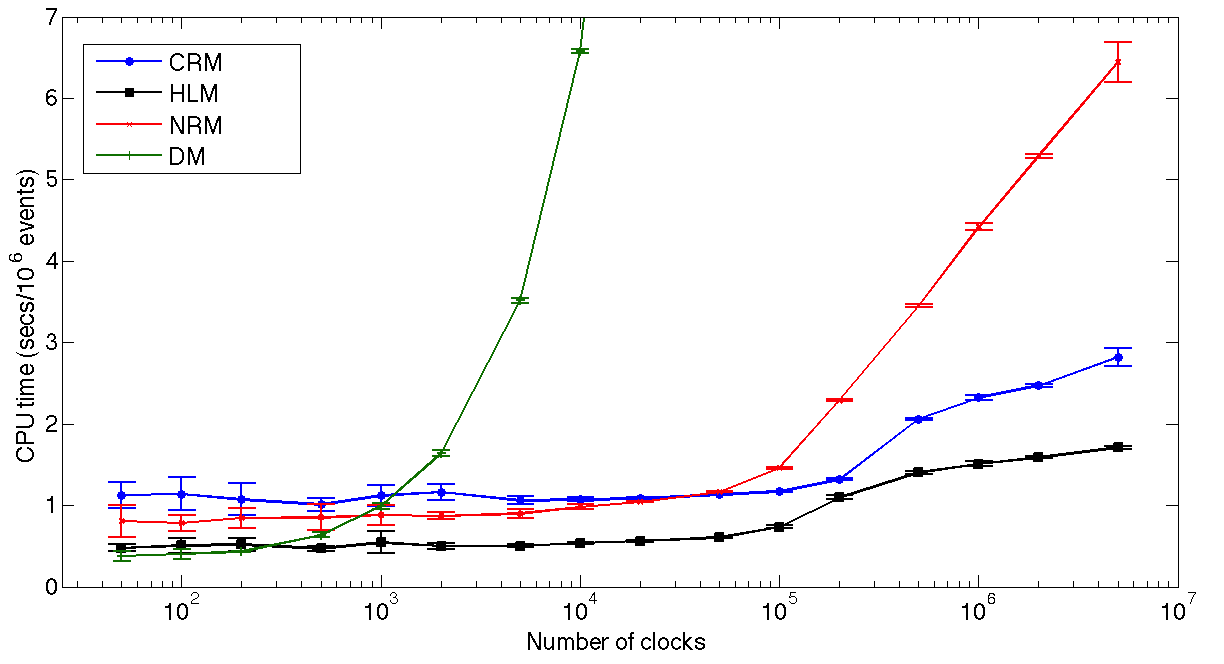}
\caption{CPU time in seconds per $10^{6}$ events vs. number of
  exponential clocks for the chemical reaction network. Blue, black, red, and green
  plots represent the mean CPU times of the CRM, the HLM, the NRM,
  and the DM over $10$ runs, respectively. The error bars indicate
  one standard deviation of the mean.}
\label{CRNspeed}
\end{figure}

\begin{figure}[h]
\centering
\includegraphics[width = \textwidth]{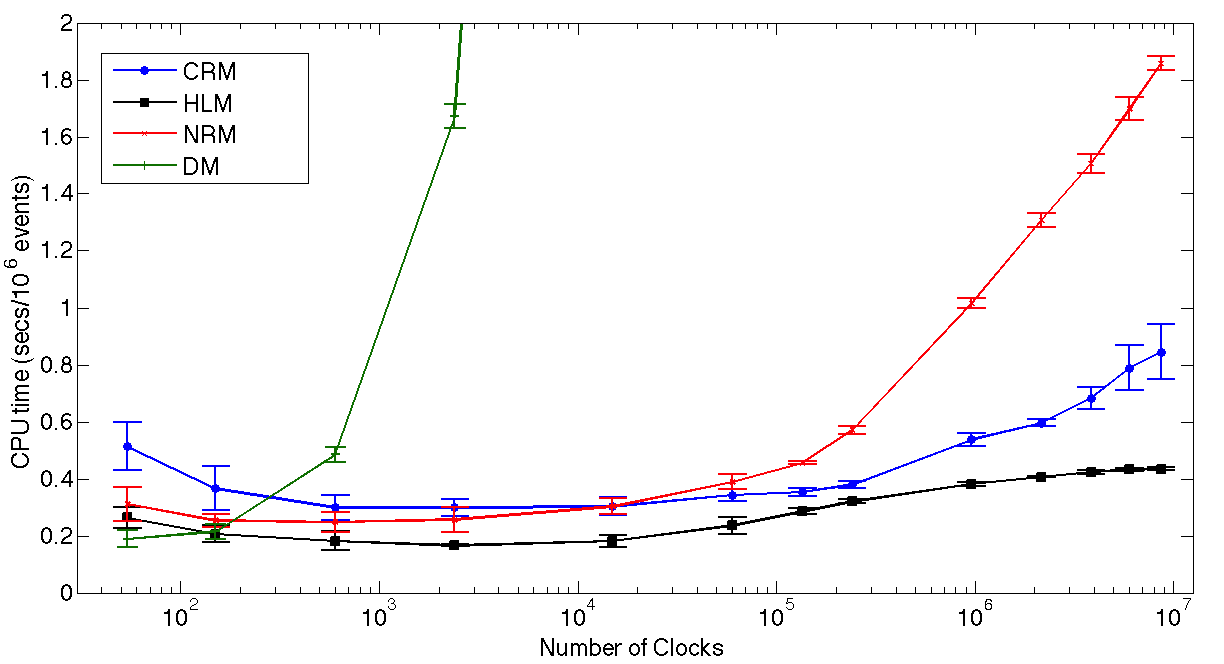}
\caption{CPU time in seconds per $10^{6}$ events vs. number of
  exponential clocks for the reaction-diffusion model. Blue, black, red, and green
  plots represent the mean CPU times of the CRM, the HLM, the NRM,
  and the DM over $10$ runs, respectively. The error bars indicate
  one standard deviation of the mean.}
\label{RDspeed}
\end{figure}

\begin{figure}[h]
\centering
\includegraphics[width = \textwidth]{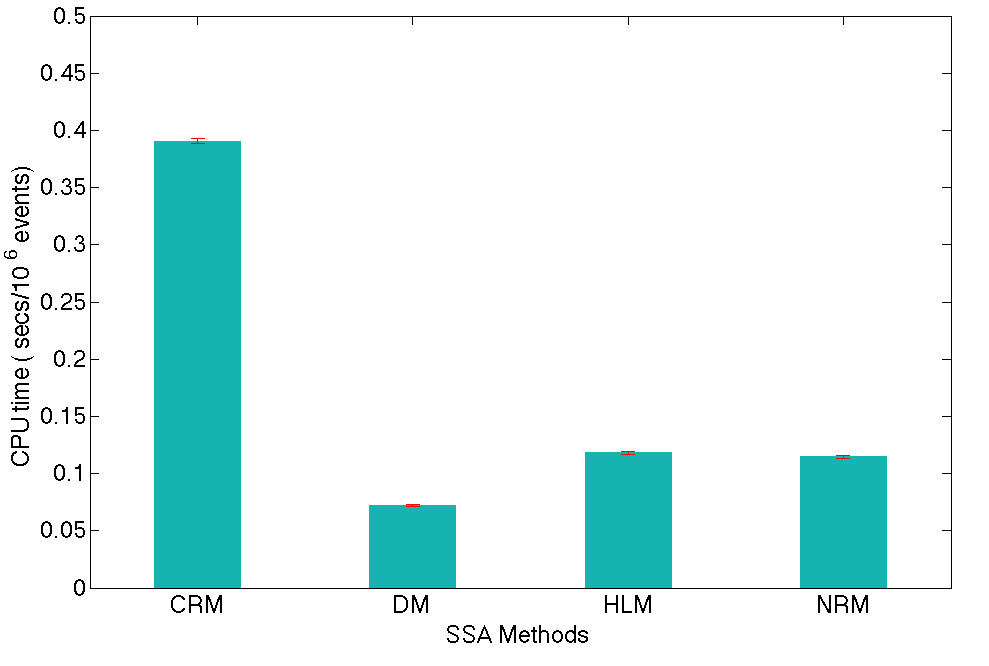}
\caption{CPU time in seconds per $10^{6}$ events vs. number of
  exponential clocks for the ``Oregonator'' model. Heights of bars
  represent the mean CPU times of the CRM, the HLM, the NRM, 
  and the DM over $10$ runs, respectively. The red error bars indicate
  one standard deviation of the mean.}
\label{Orspeed}
\end{figure}

In examples 5.1 {\bf (i) -- (iii)} , CPU times ( seconds per million events)
vs. $M$ of all four SSA methods are plotted in linear-log plots and
presented in Figure \ref{KMPspeed}, \ref{CRNspeed}, and 
\ref{RDspeed}. In these figures, each plot represents the mean CPU times of $10$
runs. The error bars indicate one standard deviation of the mean. The CPU
times ( seconds per million events) of four SSA methods over $10$ runs of the ``Oregonator''
are presented in Figure \ref{Orspeed}. The error bars also represent one
standard deviation.

As shown in these figures, for large scale models, the HLM is superior
to the other SSA methods in all tested examples. Although having the same computational
complexity, we find that the HLM is slightly faster than the CRM in all examples
and all system sizes. One possible reason is that maintaining groups
in the CRM is more expensive than maintaining buckets in the HLM. The
CRM also generates more random numbers than the HLM. For
small scale models, we find that the performance of the HLM is slower than the DM
(and sometimes the NRM) but remains competitive. 
 
{\bf (ii). Statistics of computer operations}

The  running time of an algorithm depends on many things beyond the efficiency of algorithms.
Details of implementation, the compiler, the operating
system, and the size of CPU cache can significantly change the
empirical CPU times. Therefore, it is important to study the average number of
operations per event of the HLM.   

\begin{figure}[h]
\centerline{\includegraphics[width = \linewidth]{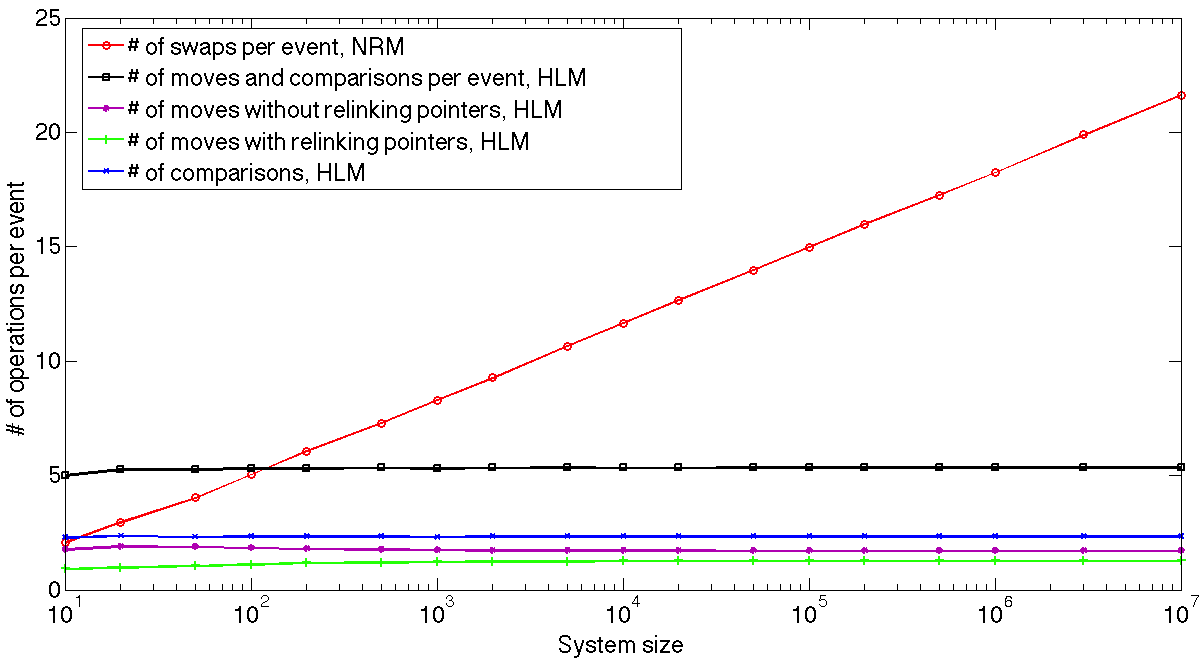}}
\caption{Number of operations per event of the generalized KMP
  model. Red: number of swaps per event of the NRM. Black: number of
  moves and comparisons per event of the HLM. Purple: number of moves
  (without relinking pointers) per event of the HLM. Green: number of
  moves (with relinking pointers) per event of the HLM. Blue: number
  of comparisons per event of the HLM. Parameters are same as in
  Section 5.2 {\bf (i)}. }
\label{stat} 
\end{figure}

We collect data of computer operations of  
  examples  5.1 {\bf (i) -- (iii) }  . The number of
comparisons when linearly searching a 
bucket and the number of moves of events between buckets  are recorded. 
Our results confirm that, in all three examples, the HLM has constant
computational cost per event. The total numbers of comparisons
plus moves per event are $\approx 5.0$ for the generalized KMP model,
$\approx 10.5$ for the reaction diffusion system, and $\approx 27.5$
for the chemical reaction network, respectively.  As an example of $O(1)$ algorithm vs. $O(\log M)$
algorithm, the statistics of computer operations of the HLM vs. the NRM for the generalized KMP
model are demonstrated in Figure \ref{stat},  in which the number of moves of
events are further broken down into moves with relinking pointers and
moves without relinking pointers (i.e. within the same bucket). The
number of operations of the NRM shown in the figure is the number of
swaps in the binary heap.  We find that the HLM is more efficient than
the NRM when $M$ is greater than $10^{2}$. 

\section{Conclusion}

In the present paper, we introduced a fast  method
 for the
  stochastic simulation algorithm (SSA), namely
the Hashing-Leaping Method (HLM), for a class of Markov jump processes arising in many
scientific fields. The common feature of these Markov jump processes is
that they are driven by many heterogeneous, state-dependent exponential clocks. The number of
exponential clocks, or the scale of the system, is denoted by $M$. As
the Markov jump process proceeds at a sequence 
of random times, the main
 strategy of the SSA is to identify the
next time of occurrence among many exponential random
times.  To do so, the HLM uses a hash-table-like data
  structure to distribute times of occurrence
covered by a certain time step with length $\tau$ into $Q$ evenly
divided buckets, updates all buckets sequentially, and leaps
forward by $\tau$.  Under assumptions (a) and (b) in Section 4.3, the average 
computational cost per event of the HLM is $O(1)$,
independent of the number of exponential clocks. 

For large scale Markov jump processes, the HLM has the desired
performance. The speed of the HLM is tested with three large-scale models: a
generalized KMP model, a chemical reaction network, and a
stochastic reaction-diffusion system. Our simulation results showed
advantages of the HLM over the DM, the NRM, and the CRM when $M$ is
greater than $\sim 10^{2}$. In addition, we verified numerically that
the average number of computer operations per event of the HLM is $O(1)$.

It is well accepted that no SSA method is unconditionally superior to
the rest. For small-scale problems, Gillespie's direct method (DM)
usually has the best efficiency. This is confirmed by our simulation
result of a chemical oscillator model called the ``Oregonator''. For this small-scale model (in which $M = 5$), the performance of the
HLM remains competitive. While possibly due to the overhead of
maintaining many groups, as the only other existing conditional $O(1)$ per-event
SSA method to the best of our knowledge, the CRM is slower than the
rest SSA methods when the scale of the Markov process is sufficiently
small.   

We do not claim that the HLM is a perfect algorithm. One drawback is
that the performance of the HLM depends on the choice of
parameters. According to our analysis in Section 4.4, the optimal
parameters depend on some constants that should be estimated
empirically. Although the performance of the HLM is
not sensitive with respect to change of parameters as long as $\tau
\sim O(1)$ and $Q \sim O(M)$, to reach the full efficiency of the
HLM, some empirical estimations or small scale CPU time tests are be
needed at the current stage. To partially solve this problem, in future, we will develop
algorithms that can adjust parameters as the simulation proceeds. 

Overall, our analysis and numerical simulations show that the HLM is a
promising SSA method, especially for large scale Markov
processes. The performance of the HLM can be potentially improved
in several aspects, from more efficient implementations to parallelizations. It is also useful to
trim the HLM for specific problems, such as multiple time
scale systems and systems with varying number of exponential clocks. Those
issues will be studied in our subsequent works.

\section*{Acknowledgement}

The authors would like to thank Lai-Sang Young, Kevin Lin, Markos
Katsoulakis, and the anonymous reviewer for many enlightening discussions
and constructive suggestions.

\bibliography{myref}
\bibliographystyle{amsplain}

\end{document}